\documentclass[11pt,leqno]{amsart}
\usepackage{graphicx}
\usepackage[english]{babel}
\usepackage{amsthm, amsmath, amssymb}
\usepackage[matrix,arrow,curve]{xy}
\newtheorem{thm}{Theorem}
\newcounter{corr}
\newtheorem{cor}[corr]{Corollary}
\newcounter{lemm}
\newtheorem{lem}[lemm]{Lemma}
\newcounter{propp}
\newtheorem{prop}[propp]{Proposition}
\newcounter{defnn}
\theoremstyle{definition}
\newtheorem{defn}[defnn]{Definition}

\theoremstyle{remark}

\numberwithin{equation}{section}
\newcommand*{\hm}[1]{#1\nobreak\discretionary{}%
{\hbox{$\mathsurround=0pt #1$}}{}}

\newcommand{\SL}{{\mathop{\mathrm{SL}}}}
\newcommand{\Sp}{{\mathop{\mathrm{Sp}}}}
\makeatletter

\renewcommand{\@listI}{%
\leftmargin=10pt \rightmargin=0pt \labelsep=5pt \itemindent=10pt \listparindent=10pt \topsep=8pt plus 2pt minus 4pt
\partopsep=2pt plus 1pt minus 1pt
\parsep=0pt plus 1pt
\itemsep=\parsep}
\begin{document}
\sloppy


\title[Parabolically connected subgroups]{Parabolically connected subgroups}%
\author{Igor V. Netay}%
\address{Department of Mechanics and Mathematics, Moscow State University, Moscow, Russia}%
\address{Independent University of Moscow, Moscow, Russia}%
\email{inetay@hse.ru}%

\thanks{The author is partially supported by AG Laboratory HSE, RF government grant, ag. 11.G34.31.0023.}%
\keywords{reductive group, parabolic subgroup, spherical subgroup, flag, Moishezon space}%

\begin{abstract}
We describe all reductive spherical subgroups of the group $\mathop{\mathrm{SL}}(n)$ which have connected intersection with any parabolic subgroup of the group $\mathop{\mathrm{SL}}(n)$. This condition guarantees that any open equivariant embedding of the corresponding homogeneous space into a Moishezon space is algebraic.
\end{abstract}
\maketitle

\section{Introduction}

Let $G$ be a reductive algebraic group over the field $\mathbb{C}$.

\defn{A closed subgroup $H\subset G$ is called {\it parabolically connected} if for any parabolic subgroup $P\subseteq G$ the intersection $P\cap H$ is connected.}\bigskip

It is useful to note that an algebraic subgroup $H$ is parabolically connected iff the intersection $B\cap H$ with any Borel subgroup $B$ is connected. Indeed, let $P\subseteq G$ be a parabolic subgroup and $B\hm\subseteq G$ be a Borel subgroup contained in $P$. Then $B$ is also a Borel subgroup of $P$. The connected algebraic group $P$ is a union of its Borel subgroups~\cite[ch.~8,~\S~22]{Humphreys},\, hence $H\cap P=\bigcup\limits_{B\subseteq P}(H\cap B)$. Since any intersection $H\cap B$ is connected and contains the identity element, we obtain that $H\cap P$ is connected.

Since any subgroup of a unipotent group is connected, we see that any unipotent subgroup $H\subset G$ is parabolically connected. It was proved by Hausen~\cite[Thm.~3]{Hausen} that for any reductive group $H$ the diagonal $\Delta H=\{(h,h)\::\:h\in H\}$ is parabolically connected as a subgroup of $G\hm=H\times H$.

Recall that an algebraic subgroup $H\subset G$ is said to be a {\it spherical subgroup} if the induced action of some Borel subgroup $B$ in $G$ on the homogeneous space $G/H$ has an open orbit. The main result of this paper is the classification of parabolically connected reductive spherical subgroups in the group $\mathop{\mathrm{SL}}(n)$. Our goal is to choose parabolically connected subgroups in the list of all spherical subgroups~\cite{Kramer}. Denote by $\mathop{\mathrm{S}}(\mathop{\mathrm{GL}}(m)\times\mathop{\mathrm{GL}}(n))$ the subgroup of $\mathop{\mathrm{SL}}(m+n)$ which consists of all block matrices with blocks of sizes $m$ and $n$. Denote by $\mathrm{T}^1$ the one-dimensional algebraic torus.

\thm{The subgroups $$\mathop{\mathrm{SL(m)}}\times\mathop{\mathrm{SL}}(n)\subset\mathop{\mathrm{SL}}(m+n)\text{ for all }m\text{ и }n,$$
$$\mathop{\mathrm{S}}(\mathop{\mathrm{GL}}(m)\times\mathop{\mathrm{GL}}(n))\subset\mathop{\mathrm{SL}}(m+n)\text{ for }m\ne n,$$
$$\mathop{\mathrm{Sp}}(2n)\subset\mathop{\mathrm{SL}}(2n),\;\mathop{\mathrm{Sp}}(2n)\subset\mathop{\mathrm{SL}}(2n+1)\text{ и }\mathop{\mathrm{Sp}}(2n)\times\mathop{\mathrm{T}}\nolimits^1\subset\mathop{\mathrm{SL}}(2n+1)$$ are parabolically connected.
At the same time the subgroups
$$\mathop{\mathrm{SO}}(n)\subset\mathop{\mathrm{SL}}(n)\text{ и
}\mathop{\mathrm{S}}(\mathop{\mathrm{GL}}(n)\times\mathop{\mathrm{GL}}(n))\subset\mathop{\mathrm{SL}}(2n)$$ are not parabolically
connected.}\label{main}\bigskip

\rm The notion of parabolically connected subgroup becomes important in complex analysis. Let $X$ be an analytic compact connected complex variety. Denote by $\mathcal{M}(X)$ the field of global meromorphic functions on $X$. It was shown that the transcendence degree of $\mathcal{M}(X)$ is less than or equal to the dimension of $X$. Those varieties for which the equality holds are called Moishezon varieties. It is known that the connected component of the identity in the group of automorphisms $\mathop{\mathrm{Aut}}^\circ(X)\subset\mathop{\mathrm{Aut}}(X)$ of a Moishezon space $X$ has a natural structure of affine algebraic group. An action of a connected reductive group $G$ on a Moishezon space $X$ is called {\it algebraic} if the corresponding homomorphism $G\to\mathop{\mathrm{Aut}}^\circ(X)$ is a homomorphism of algebraic groups. It is natural to conjecture that if the group $\mathop{\mathrm{Aut}}^\circ(X)$ is "sufficiently large", then the space $X$ is an algebraic manifold. The first result is this direction was obtained in the paper~\cite{DLuna} by D.~Luna.

\thm~\cite[Thm.~1]{DLuna} {Let $X$ be a Moishezon space equipped with an action of an algebraic torus $\mathrm{T}$ with an open orbit. Then $X$ is an algebraic $\mathrm{T}$-manifold.}\label{Luna}\bigskip

\rm The following result by J.~Hausen generalizes Theorem~\ref{Luna}.

\thm~\cite[Thm.~2]{Hausen} {Let $X$ be a compact Moishezon space, and suppose that a reductive group $G$ acts on $X$ algebraically. If for some Borel subgroup $B$ of $G$ the orbit $B\cdot x_0$ is dense in $X$ and each closed $G$-orbit contains a point $x$ such that there is a parabolic subgroup $Q\subset G$ opposite to $G_x$ with $B\subset Q$ and $G_{x_0}\cap Q$ is connected, then $X$ is a complex algebraic $G$-variety.}\bigskip

\cor{Let $H\subset G$ be a spherical parabolically connected subgroup and let $G/H\hm\to X$ be an open equivariant embedding into Moishezon space
$X$ with an action of $G$. Then $X$ is an algebraic $G$-variety.}\bigskip

\rm In many cases this corollary gives affirmative answer to the question stated in the paper~\cite{DLuna}: is it true that any Moishezon space equipped with a locally transitive action of a semisimple simply connected algebraic group $G$ is algebraic $G$-manifold if stabilizer of a point in dense orbit is connected? Example of nonalgebraic $\mathrm{PSL}(2)$-quasihomogeneous Moishezon space is constructed in the paper~\cite{Luna}. It would be interesting to find out if homogenious spaces $\mathop{\mathrm{SL}}(n)/H$ possess open equivariant embeddings into nonalgebraic Moishezon spaces where $H$ is one of two non parabolically connected reductive spherical subgroups in $\mathop{\mathrm{SL}}(n)$.

The final version is to be published in Mat.~Sb.~\cite{Me}.

The author is grateful to scientific adviser I.V.~Arjantsev for constant attention to this work.

\section{Lemmas on compatible bases}

Some statements about existence of appropriate bases will be useful to treat intersection of the subgroup $H$ and Borel subgroups. These results may be interesting themselves.

Denote the full flag $\{0\subset V_1\subset\dots\subset V_n=V\}$ in a vector space $V$ by the symbol $V_\bullet$.

\defn{We say that the basis $\{e_1,\dots,e_n\}$ of the space $V$ is {\it compatible with the flag} $V_\bullet$, if any subspace $V_i$ is spanned by some subset of this basis.}

\defn{We say that the basis $\{e_1,\dots,e_n\}$ of the space $V$ is {\it compatible with the decomposition} $V=U\oplus W$, if any vector $e_i$ lays in one of the spaces $U$ and $W$.}

\defn{The {\it hyperbolic basis} w.~r.~t. the skew-symmetric form $\omega$ is the bases $\{e_1,\dots,e_n\}$ such that for all $e_i$ the equation $\omega(e_i,\cdot)\equiv0$ holds or there exists an unique vector $e_j$ such that $\omega(e_i,e_j)\hm=\pm1$.}

\defn{Let $V_\bullet$ be a flag in the space $V$ and $W\hm\subset V$ be a subspace. Then the {\it quotient flag} $V_\bullet/W$ is the flag in $V/W$ that consists of quotient spaces $V_i/(V_i\cap W)$.}

\lem{Let $V=U\oplus W$ and $V_\bullet$ be a full flag in the space $V$. Then there exist bases $\{e_1,\dots,e_n\}$ and $\{v_1,\dots,v_n\}$ of the space $V$ such that the basis $\{e_1,\dots,e_n\}$ is compatible with the decomposition $V=U\oplus W$, the basis $\{v_1,\dots,v_n\}$ is compatible with the flag $V_\bullet$ and each vector $v_i$ equals some vector $e_l$ or the sum $e_j+e_k$ of some $e_j\in U$ and $e_k\in W$.}\label{glglbasis}

\begin{proof}

Let us construct the bases by induction. On the $k$-th step we construct the basis $\{v_1,\dots,v_k\}$ of the space $V_k$ and the basis $\{e_1,\dots,e_l\}$ of the space $\mathop{\mathrm{pr}}_U(V_k)\oplus\mathop{\mathrm{pr}}_W(V_k)$ which is compatible with the decomposition, where $\mathop{\mathrm{pr}}_U$ and $\mathop{\mathrm{pr}}_W$ are the projections onto $U$ and $W$ along $W$ and $U$ respectively. Suppose that $k$ steps of the construction are done. Let us do the $(k+1)$-st step. Note that for $i=1,\dots,n$ the equation $\dim(V_i)\hm=\dim(V_i\cap U)\hm+\dim(\mathop\mathrm{pr}_W(V_i))=\dim(V_i\cap W)+\dim(\mathop{\mathrm{pr}}_U(V_i))$ holds. Exactly one of the following four cases holds:

\begin{enumerate}
\item $$\dim(\mathop{\mathrm{pr}}\nolimits_U(V_{k+1}))=\dim(\mathop{\mathrm{pr}}\nolimits_U(V_k))+1,$$ $$\dim(\mathop{\mathrm{pr}}\nolimits_W(V_{k+1}))=\dim(\mathop{\mathrm{pr}}\nolimits_W(V_k))+1.$$ Then there exists a vector $v\in V_{k+1}$ such that $\mathop{\mathrm{pr}}_U(v)\notin\mathop{\mathrm{pr}}_U(V_k)$, $\mathop{\mathrm{pr}}_W(v)\notin\mathop{\mathrm{pr}}(V_k)$. Let us determine $v_{k+1}=v$, $e_{l+1}=\mathop{\mathrm{pr}}_U(v)$, $e_{l+2}=\mathop{\mathrm{pr}}_W(v)$.

\item $$\dim(\mathop{\mathrm{pr}}\nolimits_U(V_{k+1}))=\dim(\mathop{\mathrm{pr}}\nolimits_U(V_k))+1,$$ $$\dim(U\cap V_{k+1})=\dim(U\cap V_k)+1.$$ Take any $v_{k+1}=e_{k+1}\in U\cap(V_{k+1}\setminus V_k)$.

\item $$\dim(\mathop{\mathrm{pr}}\nolimits_W(V_{k+1}))=\dim(\mathop{\mathrm{pr}}\nolimits_W(V_k))+1,$$ $$\dim(W\cap V_{k+1})=\dim(W\cap V_k)+1.$$
This case is similar to the previous.

\item $$\dim(U\cap V_{k+1})=\dim(U\cap V_k)+1,$$ $$\dim(W\cap V_{k+1})=\dim(W\cap V_k)+1.$$

Since $V_{k+1}\cap U\subset\mathop{\mathrm{pr}}_U(V_{k+1})=\mathop{\mathrm{pr}}_U(V_k)$, there exists $u\hm\in V_k$ such that $\mathop{\mathrm{pr}}_U(u)\hm\in
V_{k+1}\setminus V_k$. Analogously, there exists $w\in V_k$: $\mathop{\mathrm{pr}}_W(w)\hm\in V_{k+1}\setminus V_k$.

In the case $\mathop{\mathrm{pr}}_U(w)\in V_{k+1}\setminus V_k$, let us determine $v=w$.

In the case $\mathop{\mathrm{pr}}_W(u)\in V_{k+1}\setminus V_k$, let us determine $v=u$.

Otherwise we say that $v=u+w$. This implies that $v\in V_k$, $\mathop{\mathrm{pr}}_U(v)\in V_{k+1}\setminus V_k$, $\mathop{\mathrm{pr}}_W(v)\in V_{k+1}\setminus V_k$. The vector $v$ can be considered as a linear combination of the basis $\{v_1,\dots,v_k\}$: $v=\sum\limits_{i=1}^k\alpha_kv_k$. Let us denote \begin{gather*} v'=\sum\limits_{i=1,\dots,k\atop\scriptstyle v_i\notin U\cup W}\alpha_kv_k. \intertext{Then we have} \mathop{\mathrm{pr}}\nolimits_U(v'-v)\hm=\sum\limits_{i=1,\dots,k\atop v_i\in U}\alpha_iv_i\in V_k,\\ \mathop{\mathrm{pr}}\nolimits_W(v'-v)=\sum\limits_{i=1,\dots,k\atop v_i\in W}\alpha_iv_i\in V_k. \intertext{Hence, }\mathop{\mathrm{pr}}\nolimits_U(v')\hm\in V_{k+1}\setminus V_k,\quad \mathop{\mathrm{pr}}\nolimits_W(v')\in V_{k+1}\setminus V_k,\\ v'=\sum\limits_{i=1,\dots,k\atop v_i\notin U\cup W}\alpha_iv_i=\sum\limits_{i=1,\dots,k_0\atop v_i\notin U\cup W}\alpha_iv_i,\quad k_0\leqslant k,\; \alpha_{k_0}\ne0.\end{gather*} Since $v_{k_0}\notin U\cup W$, we have $v_{k_0}=e_s+e_t$ for some $e_s\in U$, $e_t\in W$ by the construction. Let us substitute $v_{k_0}$, $e_s$, $e_t$ for $v'$, $\mathop{\mathrm{pr}}_U(v')$, $\mathop{\mathrm{pr}}_W(v')$. This substitution is compatible with the flag and the decomposition because $v_{k_0}\in V_{k_0}\setminus V_{k_0-1}$. This shows that required bases for $V_{k+1}$ and $\mathop{\mathrm{pr}_U}(V_{k+1})\hm\oplus\mathop{\mathrm{pr}_W}(V_{k+1})$ are constructed. \end{enumerate} \end{proof}

\lem{Let $V$ be a $2n$-dimensional vector space with the full flag $V_\bullet$ and $\omega$ be a nondegenerate skew-symmetric form in the space $V$. Then there exists a basis $\{e_1,\dots,e_{2n}\}$ in the vector space $V$ that is compatible with the flag $V_\bullet$ and hyperbolic w.~r.~t.~$\omega$.}

\begin{proof}

The proof is by induction on $n$. The basis of induction is the case $n=1$. Suppose that $V_1=\langle e_1\rangle$. Since the form $\omega$ is nondegenerate there exists a vector $e_2\in V_2=V$ such that $\omega(e_1,e_2)=1$. Antisymmetry of the form $\omega$ implies that $e_2\notin V_1$.

Suppose that inductive hypothesis is proved for $m<n$. Take any vector $e_1\in V_1\setminus\{0\}$. Determine $k=\min\{l:\omega(e_1,{\cdot})|_{V_l}\not\equiv0\}$. Choose a vector $v_k\in V_k$ such that $\omega(v_1,v_k)\hm=1$. The intersections of the subspaces of the flag $V_\bullet$ except $V_1$ and $V_k$ with the space $\langle e_1,e_k\rangle^\bot$ form a flag denoted by $V_\bullet'$. In the space $\langle e_1,e_k\rangle^\bot$ there is a basis $\{e_2,\dots,e_{k-1},e_{k+1},\dots,e_{2n}\}$ compatible with the flag $V_\bullet'$ and hyperbolic w.~r.~t.~the restriction of the form $\omega$ to $\langle e_1,e_k\rangle^\bot$. To conclude the proof, it remains to note that the basis $\{e_1,\dots,e_{2n}\}$ is as required. \end{proof}

\lem{Let $V$ be a $(2n+1)$-dimensional vector space, $U\hm\subset V$ be a hyperplane, $V_\bullet$ be a full flag in $V$ and $\omega$ be a skew-symmetric form with nondegenerate restriction to $U$. Then in $V$ there is a basis $\{e_1,\dots,e_{2n+1}\}$ such that $e_{2n+1}\in\ker(\omega)$, $e_1,\dots,e_{2n}\in U$ and each subspace $V_i$ is spanned by some vectors $e_i$, $i=1,\dots,2n$ and $e_i+e_{2n+1}$, $i=1,\dots,2n$ and the basis $\{e_1,\dots,e_{2n+1}\}$ is hyperbolic with respect to the form $\omega$. Denote $v_i=e_i$ or $v_i=e_i+e_{2n+1}$ in these cases.}\label{sp2nsl2n+1basis}

\begin{proof}

Since the space $V$ is odd-dimensional, the form $\omega$ is degenerate. Suppose that $\ker(\omega)\hm=\langle e_{2n+1}\rangle$. The restriction $\omega|_U$ is nondegenerate, so we have $e_{2n+1}\notin U$, this means that $V=U\oplus\langle e_{2n+1}\rangle$. Let elements of the full flag $U_\bullet$ be images of projection $\mathop{\mathrm{pr}}\colon V\to U$ along $\langle e_{2n+1}\rangle$ for elements of the flag $V_\bullet$. The application of the previous lemma yields existance of a basis $\{u_1,\dots,u_{2n}\}$ that is compatible with the flag $U_\bullet$ and hyperbolic w.~r.~t.~$\omega|_U$. Let vectors $v_1,\dots,v_{2n}$ be preimages of $u_1,\dots,u_{2n}$ under the projection $\mathop{\mathrm{pr}}$ such that $v_i\in V_i$. The basis $v_1,\dots,v_{2n},e_{2n+1}$ is compatible with the flag $V_\bullet$ and hyperbolic w.~r.~t.~$\omega$, however it may be not compatible with decomposition $V=U\oplus\langle e_{2n+1}\rangle$. Since $\ker(\mathop{\mathrm{pr}})=\langle e_{2n+1}\rangle$ for any $v$, we have $v-\mathop{\mathrm{pr}}(v)\in\langle e_{2n+1}\rangle$. Suppose that $v_i-\mathop{\mathrm{pr}}(v_i)=\alpha_i e_{2n+1}$. Consider indexes $i<j$ such that $\omega(e_i,e_j)=1$. One of the following four cases holds:

\begin{itemize}

\item $\alpha_i=0$, $\alpha_j=0$. Let us determine $v_i'=v_i$, $v_j'=v_j$.

\item $\alpha_i=0$, $\alpha_j\ne0$. Let us determine $v_i'=\alpha_jv_i$, $v_j'=\alpha_j^{-1}v_j$.

\item $\alpha_i\ne0$, $\alpha_j=0$. Let us determine $v_i'=\alpha_i^{-1}v_i$, $v_j'=\alpha_iv_j$.

\item $\alpha_i\ne0$, $\alpha_j\ne0$. Let us determine $v_i'=\alpha_i^{-1}v_i$, $v_j'=\alpha_iv_j-\alpha_jv_i$.

\end{itemize}

Now for $i=1,\dots,2n$ determine $e_i=\mathop{\mathrm{pr}}(v_i')\in U$. Thus the basis $\{e_1,\dots,e_{2n+1}\}$ is required. \end{proof}

\section{Cases $\mathop{\mathrm{SL}}(n)\times\mathop{\mathrm{SL}}(m)\subset\mathop{\mathrm{SL}}(m+n)$ and $\mathop{\mathrm{S}}(\mathop{\mathrm{GL}}(m)\times\mathop{\mathrm{GL}}(n))\subset\mathop{\mathrm{SL}}(m+n)$}

\rm \prop{The subgroup $\mathop{\mathrm{GL}}(m)\times\mathop{\mathrm{GL}}(n)\subset\mathop{\mathrm{GL}}(m+n)$ is parabolically connected.} \label{glgl}

\begin{proof}

The proof is by induction on $(m,n)$ by assumption that $(m',n')\leqslant(m,n)$ iff $m'\leqslant m$ and $n'\leqslant n$. The inductive basis for $m=0$ or $n=0$ is equivalent that a Borel subgroup $B\subset\mathrm{GL}$ is connected.

The following obvious remark is needed for the sequel. If $\varphi\colon G_1\to G_2$ is a surjective homomorphism of algebraic groups, the group $G_2$ is connected and the group $\ker(\varphi)$ lies in the connected component of the identity of the group $G_1^\circ$, then the group $G_1$ is connected.

Fix a decomposition $V=U\oplus W$, $\dim(U)=m$, $\dim(W)=n$ and a full flag $V_\bullet$ in the vector space $V$. Denote by $G$ the group $\mathop{\mathrm{GL}}(U)\times\mathop{\mathrm{GL}}(W)$ and by $H$ the group $G\cap\mathop{\mathrm{Stab}}(V_\bullet)$.

Choose a bases $\{e_1,\dots,e_{m+n}\}$ and $\{v_1,\dots,v_{m+n}\}$ by Lemma~\ref{glglbasis}. Renumber elements $\{e_1,\dots,e_{m+n}\}$ to satisfy $e_1,\dots,e_m\hm\in U$, $e_{m+1},\dots,e_{m+n}\in W$ and to save the ordering of elements in $U$ and the ordering of elements in $V$.

Suppose that $v_1=e_1\in U$ (the case $v_1=e_{m+1}\in W$ is similar), $U'=U/\langle e_1\rangle$ and $V_\bullet'=V_\bullet/\langle e_1\rangle$. Determine a projection $\varphi\colon H\to\mathop{\mathrm{GL}}(U'\oplus W)$. The kernel $\ker(\varphi)$ consists of matrices of the form $$\left(\begin{tabular}{c|c}$\begin{matrix}*&*\\0&E\end{matrix}$&$0$\\\hline$0$&$E$\end{tabular}\right)$$ w.~r.~t.~basis $\{e_1,\dots,e_{m+n}\}$, so it is connected. Then connectivity of $H$ follows from connectivity of the image $(\mathop{\mathrm{GL}}(U')\times\mathop{\mathrm{GL}}(W))\cap\mathop{\mathrm{Stab}}(V_\bullet')$ that is connected by the inductive assumption for $(m-1,n)$.

Suppose that $v_1=e_1+e_{m+1}$, $e_1\in U$, $e_{m+1}\in W$, $U'=U/\langle e_1\rangle$, $W'\hm=W/\langle e_{m+1}\rangle$, the projection $\varphi:H\to\mathop{\mathrm{GL}}(U'\oplus W')$ and $V_\bullet'=V_\bullet/\langle e_1,e_{m+1}\rangle$. The kernel consists of matrices of the form $$\left(\begin{tabular}{c|c}$\begin{matrix}\lambda&*\\0&E\end{matrix}$&{\LARGE$0$}\\\hline{\LARGE$0$}&$\begin{matrix}\lambda&*\\0&E\end{matrix}$\end{tabular}\right)$$ and is connected. The image equals $(\mathop{\mathrm{GL}}(U')\times\mathop{\mathrm{GL}}(W'))\cap\mathop{\mathrm{Stab}}(V_\bullet')$. Thus connectivity of $H$ follows from connectivity of the image, i.~e.~by inductive assumption for $(m\hm-1,n-1)$. \end{proof}

\prop{The subgroup $\mathop{\mathrm{SL}}(m)\times\mathop{\mathrm{GL}}(n)\subset\mathop{\mathrm{GL}}(m+n)$ is parabolically connected.} \label{slgl}

\begin{proof}

Let us prove this proposition as above in the following terms: $G=\mathop{\mathrm{SL}}(U)\hm\times\mathop{\mathrm{GL}}(W)\hm\subset\mathop{\mathrm{GL}}(U\oplus W)$, $H=G\cap\mathop{\mathrm{Stab}}(V_\bullet)$, where $(m,n)=(\dim(U),\dim(W))$, $V_\bullet$ is a full flag in the space $V=U\oplus W$. The proof is by induction on $(m,n)$ with the same ordering. The inductive basis for $m=0$ or $n=0$ is that Borel subgroups $\mathop{\mathrm{SL}}$ и $\mathop{\mathrm{GL}}$ are connected. Choose the bases $\{e_1,\dots,e_{m+n}\}$ и $\{v_1,\dots,v_{m+n}\}$and renumber in the same way.

Suppose that $v_1=e_1\in U$, $U'=U/\langle e_1\rangle$, $\varphi\colon H\to\mathop{\mathrm{GL}}(U'\oplus W)$, и $V_\bullet'=V/\langle e_1\rangle$. The kernel $\ker(\varphi)$ consists of matrices of the form $$\left(\begin{tabular}{c|c}$\begin{matrix}1&*\\0&E\end{matrix}$&$0$\\\hline$0$&$E$\end{tabular}\right).$$ The kernel is connected. The image equals $(\mathop{\mathrm{GL}}(U')\times\mathop{\mathrm{GL}}(W))\cap\mathop{\mathrm{Stab}}(V_\bullet')$ and is connected by Theorem~\ref{glgl}.

Suppose that $v_1=e_{m+1}\in W$, $W'=W/\langle e_{m+1}\rangle$, $\varphi\colon H\to\mathop{\mathrm{GL}}(U\oplus W')$, и $V_\bullet'=V_\bullet/\langle e_{m+1}\rangle$. The kernel $\ker(\varphi)$ consists of matrices of the form $$\left(\begin{tabular}{c|c}$E$&$0$\\\hline$0$&$\begin{matrix}*&*\\0&E\end{matrix}$\end{tabular}\right),$$ so it is connected. The image equals $(\mathop{\mathrm{SL}}(U)\times\mathop{\mathrm{GL}}(W'))\cap\mathop{\mathrm{Stab}}(V_\bullet')$ and is connected by the inductive assumption for $(m,n-1)$.

Suppose that $v_1=e_1+e_{m+1}$, $e_1\in U$, $e_{m+1}\in W$, $U'=U/\langle e_1\rangle$, $W'\hm=W/\langle e_{m+1}\rangle$, $\varphi\colon H\to\mathop{\mathrm{GL}}(U'\oplus W')$ and $V_\bullet'=V_\bullet/\langle e_1,e_{m+1}\rangle$. The kernel $\ker(\varphi)$ consists of matrices of the form $$\left(\begin{tabular}{c|c}$\begin{matrix}1&*\\0&E\end{matrix}$&{\LARGE$0$}\\\hline{\LARGE$0$}&$\begin{matrix}1&*\\0&E\end{matrix}$\end{tabular}\right)$$ and is connected. The image equals $(\mathop{\mathrm{GL}}(U')\times\mathop{\mathrm{GL}}(W'))\cap\mathop{\mathrm{Stab}}(V_\bullet')$ and is connected by proposition~\ref{glgl}. Therefore the group $H$ is connected. \end{proof}

\prop{The subgroup $$\{(A,B)\in\mathop{\mathrm{GL}}(m)\times\mathop{\mathrm{GL}}(n):\det(A)=\det(B)\}\hm\subset\mathop{\mathrm{GL}}(m+n)$$ is parabolically connected.} \label{glgldet}

\begin{proof}

We prove this in the same way in the following terms: $G=\{(A,B)\hm\in\mathop{\mathrm{GL}}(U)\times\mathop{\mathrm{GL}}(W):\det(A)=\det(B)\}\subset\mathop{\mathrm{GL}}(U\oplus W)$, $H\hm=G\cap\mathop{\mathrm{Stab}}(V_\bullet)$ where $(m,n)=(\dim(U),\dim(W))$, $V_\bullet$ is a full flag in the space $V=U\oplus W$. The proof is by induction on $(m,n)$ with the same ordering. The inductive basis for $m=0$ or $n=0$ is that Borel subgroups in the group $\mathop{\mathrm{SL}}$ are connected. Choose and renumber bases $\{e_1,\dots,e_{m+n}\}$ and $\{v_1,\dots,v_{m+n}\}$ in the same way.

Suppose that $v_1=e_1\in U$ (the case $v_1=e_{2n+1}\in W$ is analogous), $\varphi\colon H\hm\to\mathop{\mathrm{GL}}(U'\oplus W)$, $U'=U/\langle e_1\rangle$,
$V_\bullet'=V_\bullet/\langle e_1\rangle$. The kernel $\ker(\varphi)$ consists of elements of the form $$\left(\begin{tabular}{c|c}$\begin{matrix}1&*\\0&E\end{matrix}$&$0$\\\hline$0$&$E$\end{tabular}\right).$$ This implies that the kernel is connected. The image $\mathop{\mathrm{Im}}(\varphi)$ equals $(\mathop{\mathrm{GL}}(U')\times\mathop{\mathrm{GL}}(W))\hm\cap\mathop{\mathrm{Stab}}(V_\bullet')$ and is connected by proposition~\ref{glgl}.

Suppose that $v_1=e_1+e_{m+1}$, $e_1\in U$, $e_{m+1}\in W$, $\varphi\colon H\to\mathop{\mathrm{GL}}(U'\oplus W')$, $U'\hm=U/\langle e_1\rangle$, $W'=W/\langle
e_{m+1}\rangle$ $V_\bullet'=V_\bullet/\langle e_1,e_{m+1}\rangle$. The kernel $\ker(\varphi)$ consists of matrices of the form $$\left(\begin{tabular}{c|c}$\begin{matrix}\lambda&*\\0&E\end{matrix}$&{\LARGE$0$}\\\hline{\LARGE$0$}&$\begin{matrix}\lambda&*\\0&E\end{matrix}$\end{tabular}\right).$$ This yields it is connected. The image $\mathop{\mathrm{Im}}(\varphi)$ equals $(\{(A,B)\in\mathop{\mathrm{GL}}(U')\hm\times\mathop{\mathrm{GL}}(W'):\det(A)=\det(B)\})\cap\mathop{\mathrm{Stab}}(V_\bullet')$ and is connected by the inductive assumption for~$(m\hm-1,n-1)$.

Therefore the group $H$ is connected. \end{proof}

\prop{The subgroup $\mathop{\mathrm{SL}}(m)\times\mathop{\mathrm{SL}}(n)\subset\mathop{\mathrm{SL}}(m+n)$ is parabolically connected.} \label{slsl}

\begin{proof}

Let us prove that the subgroup $G=\mathop{\mathrm{SL}}(m)\times\mathop{\mathrm{SL}}(n)\hm\subset\mathop{\mathrm{GL}}(m+n)$ is parabolically connected. Let $B\subset\mathop{\mathrm{GL}}(m+n)$ be a Borel subgroup. Then the group $B'=B\hm\cap\mathop{\mathrm{SL}}(m+n)$ is a Borel subgroup for $\mathop{\mathrm{SL}}(m+n)$, и $G\cap B=G\cap B'$, because $G\subset\mathop{\mathrm{SL}}(m+n)$. Hence parabolic connectivity of the group $\mathop{\mathrm{SL}}(n)\times\mathop{\mathrm{SL}}(m)$ as a subgroup of $\mathop{\mathrm{SL}}(m+n)$ is equivalent to parabolic connectivity as a subgroup of $\mathop{\mathrm{GL}}(m+n)$.

The further proof and terms are similar to the previous: $G=\mathop{\mathrm{SL}}(U)\hm\times\mathop{\mathrm{SL}}(W)$, $H\hm=G\cap \mathop{\mathrm{Stab}}(V_\bullet)$ where $(m,n)=(\dim(U),\dim(W))$, $V_\bullet$ is a full flag in the space $V\hm=U\oplus W$. The proof is by induction on $(m,n)$ with the same ordering. The inductive basis for $m=0$ or $n=0$ is that Borel subgroups in the group $\mathop{\mathrm{SL}}$ are connected. Choose and and renumber bases $\{e_1,\dots,e_{m+n}\}$ and $\{v_1,\dots,v_{m+n}\}$ as above.

Suppose that $v_1=e_1\in U$, $U'=U/\langle e_1\rangle$, $\varphi\colon H\to\mathop{\mathrm{GL}}(U'\oplus W)$ and $V_\bullet'=V_\bullet/\langle e_1\rangle$. One of three following cases holds:

\begin{itemize}

\item $\forall i = 2,\dots,m+n\quad\exists j\in\{2,\dots,m+n\}\quad\exists k: v_k=e_i+e_j$. Denote by $s$ a permutation of the set $\{2,\dots,m+n\}$ such that for any $i$, $j$, $k$ satisfying $e_i+e_j=v_k$ we have $s(i)=j$ and $s(j)=i$. Assume an element $g\in G$ has matrix $(a_{ij})$ with respect to the basis $\{e_1\dots,e_{m+n}\}$. Let us show that $a_{ii}\hm=a_{s(i)s(i)}$ for $i=2,\dots,m+n$. Fix $i$ and $v_k=e_i+e_{s(i)}$. For any element $g\in\mathop{\mathrm{Stab}}(V_\bullet)$ we have $gv_k=\lambda_kv_k+v$, $v\in V_{k-1}$. Also, $V_k=V_{k-1}\oplus\langle v_k\rangle$ and $V_{k-1}\cap\langle e_i,e_{s(i)}\rangle=0$ by the construction of the basis, so $a_{ii}=a_{s(i)s(i)}=\lambda_k$. Suppose that $\{k:v_k\notin U\cup W\}=\{i_1,\dots,i_l\}$. Then $\det(g|_U)=a_{11}\lambda_{i_1}\dots\lambda_{i_l}\hm=1$ and $\det(g|_W)=\lambda_{i_1}\dots\lambda_{i_l}$. Therefore $g\in\mathop{\mathrm{SL}}(U)\times\mathop{\mathrm{SL}}(W)$. This yields that $a_{11}=1$. Hence, the kernel consists of matrices of the form $$\left(\begin{tabular}{c|c}$\begin{matrix}1&*\\0&E\end{matrix}$&$0$\\\hline$0$&$E$\end{tabular}\right),$$ and is connected. The image equals $(\mathop{\mathrm{SL}}(U')\times\mathop{\mathrm{SL}}(W))\cap\mathop{\mathrm{Stab}}(V_\bullet')$ and is connected by the inductive hypothesis.

\item $\exists e_i=v_j\in U$, $\nexists i',j':v_{j'}=e_i+e_{i'}$. Then the group $H$ contains the one-dimensional torus $T=\mathop{\mathrm{diag}}(\lambda,1,\dots,1,\lambda^{-1},1,\dots,1)$, where $\lambda^{-1}$ equals the $i$-th coordinate. Multiplying by $t$ preimages of all elements of the group $(\mathop{\mathrm{SL}}(U')\times\mathop{\mathrm{SL}}(W))\cap\mathop{\mathrm{Stab}}(V_\bullet')$, we get preimages of all elements for the group $(\mathop{\mathrm{GL}}(U')\times\mathop{\mathrm{SL}}(W))\cap\mathop{\mathrm{Stab}}(V_\bullet')$. As before, the kernel is a unipotent group and then in particalar is connected. The image is connected by proposition~\ref{slgl}.

\item $\exists e_i=v_j\in W$: $\exists!i',j': v_{j'}=e_i+e_{i'}$. It can be assumed that $\dim(U)>1$. In the converse case, $G=\mathop{\mathrm{SL}}(W)$ and the inductive statement is that a Borel subgroup for the group $\mathop{\mathrm{SL}}$ is connected. Suppose that $e_s+e_k=v_l$. Then the group $H$ contains the one-dimensional torus $T\hm=\{\mathop{\mathrm{diag}}(\lambda,1,\dots,1,\lambda^{-1},1,\dots,1,\lambda^{-1},1,\dots,1,\lambda,1,\dots,1)\}$ where the value $\lambda$ coincides with the first and $i$-th coordinates and the value $\lambda^{-1}$ coincides with $s$-th and $k$-th coordinates. As in the previous case, the image equals $(\mathop{\mathrm{GL}}(U')\times\mathop{\mathrm{SL}}(W))\cap\mathop{\mathrm{Stab}}(V_\bullet')$. If there are not $s$, $k$, $l$ such that $e_k+e_s=v_l$, then assume $s=2$, $k>m$, $k\ne i$. The torus $T$ lies in the group $H$. The kernel and the image are the same to previous.

\end{itemize}

Suppose that $v_1=e_1+e_{m+1}$, $e_1\in U$, $e_{m+1}\in W$, $U'=U/\langle e_1\rangle$, $W'=W/\langle e_{m+1}\rangle$, $\varphi:H\to\mathop{\mathrm{GL}}(U'\oplus W')$, и $V_\bullet'=V_\bullet/\langle e_1,e{m+1}\rangle$. We may assume that $\dim(U)>1$ and $\dim(W)>1$. Otherwise the statement is that a Borel subgroup for $\mathrm{SL}$ is connected. Then we have either $\exists v_k=e_i+e_j$, $i>1$, $j>m+1$, $e_i\in U$, $e_j\in W$ or $\exists v_k=e_i\in U$ and $\exists v_l=e_j\in W$, $i>1$, $j>m+1$. In both cases the group $H$ contains one-dimensional torus $T=\{\mathop{\mathrm{diag}}(\lambda,1,\dots,1,\lambda^{-1},1,\dots,1,\lambda,1,\dots,1,\lambda^{-1})\}$ where the value $\lambda$ coincides the first and $(m+1)$-th coordinate, the value $\lambda^{-1}$ coincides $i$-th and $j$-th coordinates. Similarly, the image equals $\{(A,B)\in\mathop{\mathrm{GL}}(U')\times\mathop{\mathrm{GL}}(W'): \det(A)=\det(B)\}\cap\mathop{\mathrm{Stab}}(V_\bullet')$. Elements in the kernel of the map $\varphi$ have the form $$\left(\begin{tabular}{c|c}$\begin{matrix}1&*\\0&E\end{matrix}$&{\LARGE$0$}\\\hline{\LARGE$0$}&$\begin{matrix}1&*\\0&E\end{matrix}$\end{tabular}\right),$$ so the kernel is connected.

Since the kernel and the image are connected, the group $H$ is connected.

\end{proof}

\section{Case $\mathop{\mathrm{S}}(\mathop{\mathrm{GL}}(m)\times\mathop{\mathrm{GL}}(n))\subset\mathop{\mathrm{SL}}(m+n)$, $m\ne n$}

\rm \prop{The subgroup $\mathop{\mathrm{S}}(\mathop{\mathrm{GL}}(m)\times\mathop{\mathrm{GL}}(n))\subset\mathop{\mathrm{SL}}(m+n)$ is parabolically connected.}

\begin{proof}

The proof is by induction on $(m,n)=(\dim(U),\dim(W))$ in the following terms: $G=\mathop{\mathrm{S}}(\mathop{\mathrm{GL}}(U)\hm\times\mathop{\mathrm{GL}}(W))$, $H=G\cap\mathop{\mathrm{Stab}}(V_\bullet)$ where $V_\bullet$ is a full flag, $V=U\oplus W$. The inductive basis for $m=0$ or $n=0$ is that a Borel subgroup of the group $B\subset\mathrm{SL}$ is connected. As above, choose and renumber bases $\{e_1,\dots,e_{m+n}\}$ and $\{v_1,\dots,v_{m+n}\}$.

Suppose that $v_1=e_1\in U$ (the case $v_1=e_{m+1}\in W$ is similar), $U'=U/\langle e_e\rangle$, $\varphi\colon H\hm\to\mathop{\mathrm{GL}}(U'\oplus W)$ and $V+\bullet'=V_\bullet/\langle e_1\rangle$. The kernel $\ker(\varphi)$ consists of matrices of the form $$\left(\begin{tabular}{c|c}$\begin{matrix}1&*\\0&E\end{matrix}$&$0$\\\hline$0$&$E$\end{tabular}\right),$$ and is connected. Since for all $g'\in G'$ one can choose a matrix element $a_{11}$ of the preimage such that the determinant of element in preimage equals $1$, the image equals $G'=(\mathop{\mathrm{GL}}(U')\times\mathop{\mathrm{GL}}(W))\cap\mathop{\mathrm{Stab}}(V_\bullet')$.

Suppose that $v_1=e_1+e_{m+1}$, $U'=U/\langle e_1\rangle$, $W'=W\langle e_{m+1}\rangle$, $\varphi\colon H\to\mathop{\mathrm{GL}}(U'\oplus W')$ and $V_\bullet'=V_\bullet/\langle e_1,e_{m+1}\rangle$. The kernel $\ker(\varphi)$ consists of matrices of the form $$\left(\begin{tabular}{c|c}$\begin{matrix}1&*\\0&E\end{matrix}$&{\LARGE$0$}\\\hline{\LARGE$0$}&$\begin{matrix}1&*\\0&E\end{matrix}$\end{tabular}\right)$$ and is connected. The image equals $(\mathop{\mathrm{GL}}(U')\times\mathop{\mathrm{GL}}(W'))\cap\mathop{\mathrm{Stab}}(V_\bullet')$ and is connected by proposition~\ref{glgl}.

Therefore, connectivity of the kernel and the image implies connectivity of the group $H$.

\end{proof}

\section{Cases $\mathop{\mathrm{Sp}}(2n)\subset\mathop{\mathrm{SL}}(2n)$, $\mathop{\mathrm{Sp}}(2n)\subset\mathop{\mathrm{SL}}(2n+1)$ and $\mathop{\mathrm{Sp}}(2n)\times\mathrm{T}^1\subset\mathop{\mathrm{SL}}(2n+1)$}

\prop{The subgroup $\mathop{\mathrm{Sp}}(2n)\subset\mathop{\mathrm{SL}}(2n+1)$ is parabolically connected.}\label{sp2nsl2n+1}

\begin{proof}
Suppose that $V=U\hm\oplus W$, $\dim(U)=2n$, $\dim(W)=1$. Let $\omega$ be a skew-symmetric form on $V$ such that the restriction $\omega|_U$ is nondegenerate, $W\hm=\ker(\omega)$, $V_\bullet$ be a full flag in the space $V$ and $H\hm=\mathop{\mathrm{Sp}}(2n)\hm\cap\mathop{\mathrm{Stab}}(V_\bullet)$. By Lemma~\ref{sp2nsl2n+1basis} choose a basis $\{e_1,\dots,e_{2n+1}\}$ in the space $V$. We can write equations on  matrix elements in terms of this basis. Assume $A\hm=(a_{ij})\hm\in\mathop{\mathrm{SL}}(V)$. The condition $A\in H$ is equivalent to the conditions $A|_U\in\mathop{\mathrm{SL}}(U)$, $A^t\Omega A\hm=\Omega$, $Av_i\in V_i$, $i=1,\dots,2n+1$, where $\Omega=(\omega(e_i,e_j))$. Since $\forall(a_{ij})\hm\in\mathop{\mathrm{Sp}}(2n)\hm\subset\mathop{\mathrm{SL}}(2n+1)$, we have $a_{2n+1,2n+1}=1$, $a_{i,2n+1}\hm=a_{2n+1,i}\hm=0$, $i=1,\dots,2n$. Hence we can consider $2n\times 2n$-matrices.

Let us introduce the following notation. Denote $I=\{1,\dots,2n\}$, \begin{equation}S=\{i:\exists j\;v_j=e_i\hm+e_{2n+1}\}.\label{defS}\end{equation} Denote by $I^\flat, I^\sharp$ subsets in $I$ such that if $\omega(e_i,e_j)=1$, then $i\in I^\flat$, $j\in I^\sharp$. The conditions that the basis $\{e_i\}_{i\in I}$ is hyperbolic and the restriction of the form $\omega$ to the hyperplane $U$ is nondegenerate imply that $I^\flat\sqcup I^\sharp=I$. Let us write $i=j^\flat$, $i^\sharp=j$ for shortness. Only one of the expressions $i^\flat$ and $i^\sharp$ makes sense. Denote this expression by $\bar{i}$. The symbol $i^\sharp$ makes sense only for $i\in I^\flat$ and similarly for $j^\flat$. Now let us obtain a system of equation on elements of the group~$H$:

\begin{align}
&a_{i,2n+1}=0,\quad i=1,\dots,2n\,\label{sp1}\\  &a_{2n+1,i}=0,\quad i=1,\dots,2n\,\label{sp2}\\  &a_{2n+1,2n+1}=1\,\label{sp3}\\ &a_{i,j}=0,\,\qquad
i>j\,\label{sp4}\\ \text{\raisebox{1em}[1pt][1pt]{$\left\{\vphantom{\begin{aligned}
a_{i,2n+1}=0,\quad i=1,\dots,2n\, (1)\\
a_{2n+1,i}=0,\quad i=1,\dots,2n\, (2)\\
a_{2n+1,2n+1}=1\, (3)\\
a_{i,j}=0,\,\qquad i>j\,(4)\\
\sum_{i\in I^\flat}a_{i,l}a_{i^\sharp,m}-\sum_{i\in I^\sharp}a_{i,l}a_{i^\flat,m}=\omega(e_l,e_m),\,l,m=1\dots,2n\,(5)\\
\sum_{i\in S}a_{i,k}=
\left\{\begin{aligned}1,\quad&k\in Z\\
0,\quad&k\notin S\end{aligned}\right.\,(6)
\end{aligned}}\right.$}}
&\sum_{i\in I^\flat}a_{i,l}a_{i^\sharp,m}-\sum_{i\in I^\sharp}a_{i,l}a_{i^\flat,m}=\omega(e_l,e_m),\quad l,m=1\dots,2n\,\label{sp5}\\ &\sum_{i\in
S}a_{i,k}= \left\{\begin{aligned}1,\quad&k\in S\\ 0,\quad&k\notin S\end{aligned}\right.\,\label{sp6}
\end{align}

Let us show how to obtain these equations. The equations $(\ref{sp1})-(\ref{sp3})$ are equivalent that $A\in\mathop{\mathrm{SL}}(W)$.

Invariance of the flag $Av_i\in V_i$ implies that $Av_i\in\langle e_1,\dots,e_i,e_{2n+1}\rangle$; $Ae_{2n+1}=e_{2n+1}$, $v_i=e_i$ or $v_i=e_i+e_{2n+1}$
for $i=1,\dots,2n\hm\Rightarrow Ae_i\in\langle e_1,\dots,e_i,e_{2n+1}\rangle$. Combining this and~(\ref{sp2}), we obtain~(\ref{sp4}).

The equations~(\ref{sp5}) are equivalent that $A^t\Omega A=\Omega$.

Suppose that $k\in S$. Then we have $Av_k=A(e_k\hm+e_{2n+1})=Ae_k+e_{2n+1}\hm\in\langle\{e_i\}_{i\notin S,i\leqslant k},\{e_i\hm+e_{2n+1}\}_{i\in S,i\leqslant k}\rangle$. At the same time the equations $Ae_k=\sum\limits_{i=1}^{k}a_{i,k}e_i$, $Ae_{2n+1}\hm=e_{2n+1}$ hold. This yields that $e_{2n+1}=\sum\limits_{S\ni i\leqslant
k}a_{i,k}e_{2n+1}$,~i.~e.~$\sum\limits_{S\ni i\leqslant k}a_{i,k}\hm=1\Leftrightarrow\sum\limits_{i\in S}a_{i,k}=1$.

By the same arguments for $k\notin S$, we have $\sum\limits_{i\in S}a_{i,k}=0$.

Now let us conclude from these equations that $A\hm\in\mathop{\mathrm{SL}}(W)$, $A^t\Omega A=\Omega$, $Av_i\in V_i$ для $i=1,\dots,2n+1$. Evidently, the system of equations follows the first and the second conditions. Let us check the third.

Suppose that $k\in S$. Then $Av_k=Ae_k+e_{2n+1}=\sum\limits_{i=1}^{k}a_{i,k}e_i+e_{2n+1}=\sum\limits_{i=1}^{k}a_{i,k}e_i\hm+\sum\limits_{S\ni i\leqslant
k}a_{i,k}e_{2n+1}=\sum\limits_{i=1}^{k}a_{i,k}v_i\in V_i$. The case $k\notin S$ is similar.

We prove the proposition by induction. The subgroup $\mathop{\mathrm{SL}}(2)\subset\mathop{\mathrm{SL}}(3)$ is parabolically connected by proposition~\ref{slsl}. \rm The inductive assumption is that the subgroup $\mathop{\mathrm{Sp}}(2n-2)\hm\subset\mathop{\mathrm{SL}}(2n-1)$ is parabolically connected. Let us show that the subgroup $\mathop{\mathrm{Sp}}(2n)\subset\mathop{\mathrm{SL}}(2n\hm+1)$ is parabolically connected. Denote by $L$ the vector space $\langle e_1,e_{1^\sharp}\rangle$. Consider the homomorphism $\varphi\colon H\to\mathop{\mathrm{GL}}(V/L)$, where $H=\mathop{\mathrm{Sp}}(2n)\cap\mathop{\mathrm{Stab}}(V_\bullet)$. We shall see that its kernel is connected. Suppose that $(a_{ij})\in\ker(\varphi)$. Then the matrix $(a_{ij})$ has the form $$\begin{pmatrix}\lambda&*&\mu&*\\0&E&*&0\\0&0&\lambda^{-1}&*\\0&0&0&E\end{pmatrix}.$$

By equation~(\ref{sp5}) we obtain that elements denoted by $*$ equals $0$. Either some subspace $V_k$ contains $e_1+e_{2n+1}$ and does not contain $e_1$ or contains $e_{1^\sharp}+e_{2n+1}$ and does not contain $e_{1^\sharp}$. Then for all $(a_{ij})\in\ker(\varphi)$ we have $a_{11}=a_{1^\sharp1^\sharp}\hm=\lambda=1$. Otherwise the value $\lambda$ may be any in $\mathbb{C}^\times$. Thus the group is unipotent or equals semidirect product of one-dimensional torus and an unipotent group. Therefore the kernel $\ker(\varphi)$ is connected. Clearly, the image lies in $\mathop{\mathrm{Sp}}(2n-\nolinebreak2)\hm\cap\mathop{\mathrm{Stab}}(V_\bullet/L)$. Let us show that the image equals $\mathop{\mathrm{Sp}}(2n-\nolinebreak2)\hm\cap\mathop{\mathrm{Stab}}(V_\bullet/L)$. Suppose that $\varphi(e_i)=e_i'$. Then the set $\{e_i'\}_{i\ne1,1^\sharp}$ is a basis of $V/L$. Let $(a_{ij}')$ be a matrix of element $A'\in\mathop{\mathrm{Sp}}(2n-\nolinebreak2)\hm\cap\mathop{\mathrm{Stab}}(V_\bullet/L)$ w.~r.~t.~this basis. Let $A$ be an element of preimage with the matrix $a_{ij}$ in the basis $\{e_i\}$ , where $a_{ij}=a_{ij}'$ for $\{i,j\}\cap\{1,1^\sharp\}=\varnothing$, $a_{11}=a_{1^\sharp1^\sharp}=1$, the rest elements equals $0$. We shall see that $A\in H$. Indeed, the operator $A$ stabilizes the skew-symmetric form $\omega$ and the flag $V_\bullet$. Since $A'$ stabilizes $V_i/L$, the operator $A$ stabilizes $V_i$. Since $A$ stabilizes $e_{1}$, the statements $Ae_{2n+1}\hm\in V_i/L$ for $e_{2n+1}\in V_i/L$ and $A(e_k+e_{2n+1})\in V_i$ for $e_k+e_{2n+1}\in V_i$, $k=1$ or $1^\sharp$ are equivalent. Hence connectivity of the group $\mathop{\mathrm{Sp}}(V/L)\hm\cap\mathop{\mathrm{Stab}}(V_\bullet')$ implies connectivity of the group $H$. \end{proof}

\prop{The group $\mathop{\mathrm{Sp}}(2n)\subset\mathop{\mathrm{SL}}(2n)$ is parabolically connected.}

\begin{proof} Let us conclude this from parabolic connectivity of $\mathop{\mathrm{Sp}}(2n)\subset\mathop{\mathrm{SL}}(2n+1)$. Indeed, assume that the group $\mathop{\mathrm{Sp}}(2n)$ acts on the hyperplane $U$ in the space $V$. Then connectivity of intersection $\mathop{\mathrm{Sp}}(2n)$ with stabilizer of any full flag  $V_\bullet$ implies connectivity of intersection with stabilizer of any flag $U_\bullet\cup\{V\}$, i.~e. connectivity of intersection with any Borel subgroup in $\mathop{\mathrm{SL}}(V)$ implies connectivity of the intersection with any Borel subgroup for $\mathop{\mathrm{SL}}(U)$.
\end{proof}

\prop{The subgroup $\mathop{\mathrm{Sp}}(2n)\times\mathrm{T}^1\subset\mathop{\mathrm{SL}}(2n+1)$ is parabolically connected, where $t(\lambda)|_U=\lambda$ and $t(\lambda)e_{2n+1}=\lambda^{-2n}e_{2n+1}$ for $\{t(\lambda)=\mathrm{diag}(\lambda,\dots,\lambda,\lambda^{-2n})\}=\mathrm{T}^1$.}

\begin{proof}

Suppose that $V=U\oplus W$, $\dim(W)=1$, $\Sp(2n)\subset\SL(U)$. Let $B$ be a Borel subgroup for $\SL(V)$. Determine a homomorphism $\varphi\colon\Sp(2n)\times T^1\to\mathbb{C}^\times$, $\varphi(A)=A|_W$. By Lemma~\ref{sp2nsl2n+1basis} choose a basis $\{e_1,\dots,e_{2n+1}\}$ in the space $V$. In notation of proposition~\ref{sp2nsl2n+1} for each pair $(i,\bar\imath)$, $i\in S$ determine $t_i=\lambda^{-n}$, $t_{\bar\imath}=\lambda^{n+1}$, $t_{2n+1}=\lambda^{-n}$, where the set $S$ is determined by formula~(\ref{defS}). If $i$, $\bar\imath\notin S$, then determine $t_i$ similarly to the case $i\in S$ or $\bar\imath\in S$. Then $\mathop{\mathrm{diag}}(t_1,\dots,t_{2n+1})\hm\in(\Sp(2n)\times T^1)\cap B$. Hence the image $\varphi(\Sp(2n)\times T^1\cap B)$ equals $\mathbb{C}^\times$. The kernel equals $\Sp(2n)\cap B$ and is connected by proposition~\ref{sp2nsl2n+1}. Therefore for any Borel subgroup $B\subset\SL(V)$ the intersection $B\cap(\Sp(2n)\times T^1)$ is connected.
\end{proof}

\rm This completes the proof of Theorem~\ref{main}.

\section{Lack of parabolic connectivity}

Assume that the group $\mathop{\mathrm{SO}}(n)$, $n\geqslant2$ stabilizes the standard symmetric form. Then its intersection with the group of superdiagonal matrices is finite and therefore is not connected. This shows that the subgroup $\mathop{\mathrm{SO}}(n)\subset\mathop{\mathrm{SL}}(n)$ is not parabolically connected.

\prop{The subgroup $\mathop{\mathrm{S}}(\mathop{\mathrm{GL}}(n)\times\mathop{\mathrm{GL}}(n))\subset\mathop{\mathrm{SL}}(2n)$ is not parabolically connected.}

\begin{proof}

Let $\{e_1,\dots,e_n\}$ be a basis in the space $U$ and $\{e_{n+1},\dots,e_{2n}\}$ be a basis in the space $W$, $V_i=V_{i-1}\oplus\langle e_i+e_{n+i}\rangle$, $i=1,\dots,n$, $V_0=0$, $V_{n+i}=V_{n+i-1}\oplus\langle e_i\rangle$, $i=1,\dots,n$. Consider the group $H=\mathop{\mathrm{S}}(\mathop{\mathrm{GL}}(U)\hm\times\mathop{\mathrm{GL}}(W))\cap\mathop{\mathrm{Stab}}(V_\bullet)$. Take an arbitrary element $g\in H$ and consider its matrix $(a_{ij})$ with respect to the basis $\{e_1,\dots,e_{2n}\}$. Sice $g(e_i+e_{n+i})\in V_i$, we have $a_{ii}+a_{i,n+i}=a_{n+i,i}+a_{n+i,n+i}$, at the same time $a_{i,n+i}=a_{n+i,i}=0$ by invariance of the subspaces $U$ and $W$. Hence we have $a_{ii}=a_{n+i,n+i}=\lambda_i$ for $i=1,\dots,n$. It follows that $\lambda_1^2\dots\lambda_n^2=\det(g)=1$, i.~e.~$\lambda_1\dots\lambda_n=\pm1$. The group consists of two non-intersecting nonempty closed subsets. This means that the group $H$ is not connected.
\end{proof}

\bibliographystyle{amsplain}

\begin{thebibliography}{XX}


\bibitem{Humphreys} J.E.~Humphreys, Linear Algebraic Groups (Springer-Verlag, Berlin, 1975)

\bibitem{Hausen}J.~Hausen, Algebraicity criteria for almost homogeneous complex spaces. Arch. Math. 74 (2000), 317-320.

\bibitem{Kramer}M.~Kr\"amer, Sph\"arische untergruppen in kompakten zusammenh\"angenden Liegruppen. Compositio Math. 38:2 (1979),
129-153.

\bibitem{DLuna}D.~Luna, Toute vari\'et\'e de Moisezon presque homog\`ene sous un tore est un sch\'ema. C.~R.~Acad, Sci.~Paris, t.~314, S\'erie~{\rm
I}, p.~65-67, 1992.

\bibitem{Luna}D.~Luna, L.~Moser-Jauslin, Th.~Vust, Almost homogeneous Artin-Moishezon varieties under the action of ${\rm PSL}\sb 2(\mathbb{C})$. Topological
methods in algebraic transformation groups (New Brunswick, NJ, 1988), 107-115, Progr. Math. 80, Birkh\"auser Boston, MA, 1989.

\bibitem{Me}I.~Netay, "Parabolically connected subgroups", to appear in Mat.~Sb.

\end{thebibliography}

\end{document}